\documentclass[12pt]{amsart}
\newcommand{\E}{{\mathbb E}}
\newcommand{\Z}{{\mathbb Z}}

\newcommand{\CI}{{\mathcal I}}
\newcommand{\CJ}{{\mathcal J}}
\newcommand{\CL}{{\mathcal L}}
\newcommand{\CW}{{\mathcal W}}
\newcommand{\CX}{{\mathcal X}}
\newcommand{\CY}{{\mathcal Y}}
\newcommand{\CZ}{{\mathcal Z}}

\newcommand{\norm}[1]{\lVert #1\rVert}
\newcommand{\nnorm}[1]{\lvert\!|\!| #1|\!|\!\rvert}
\newcommand{\inv}{^{-1}}
\DeclareMathOperator{\id}{Id}

\newcommand{\one}{{\mathbf 1}}

\newtheorem{theorem}{Theorem}
\newtheorem*{theorem*}{Theorem}
\newtheorem{proposition}{Proposition}
\newtheorem{lemma}{Lemma}
\newtheorem{corollary}{Corollary}
\theoremstyle{definition}
\newtheorem*{definition}{Definition}

\theoremstyle{remark}
\newtheorem{remark}{Remark}

\begin{document}
\title[Ergodic seminorms for 
commuting transformations]{Ergodic seminorms for 
commuting transformations and applications}

\author{Bernard Host}

\address{ Universit\'e Paris-Est, Laboratoire d'analyse et de
math\'ematiques appliqu\'ees, UMR CNRS 8050, 5 bd Descartes, 77454 Marne la
Vall\'ee Cedex 2, France}
\email{bernard.host@univ-mlv.fr}

\date{\today}

\begin{abstract}
Recently, T. Tao gave a finitary proof a convergence theorem for multiple averages 
with several commuting transformations and soon later, T.~Austin gave an ergodic proof 
of the same result. 
Although we give here one more proof of the same theorem, this is not the 
main goal of this paper. Our main concern is to provide some tools 
for the case of several commuting transformations, similar to the tools 
successfully used in the case of a single transformation, 
with the idea that they will be 
useful in the solution of other problems.
\end{abstract}

\maketitle

\section{Introduction}
\subsection{Motivation and context}
Recently, T. Tao~\cite{T} proved a convergence result for several 
commuting transformations.
\begin{theorem*}[T. Tao]
Let $(X,\mu,S_1,\dots,S_d)$ be a system where $S_1,\dots,S_d$ are 
commuting measure preserving transformations.  Then, for every 
$f_1,\dots,f_d\in L^\infty(\mu)$, the averages
\begin{equation}
\label{eq:av}
\frac 1N\sum_{n=0}^{N-1} S_1^nf_1\cdot\ldots\cdot S_d^nf_d
\end{equation}
converge in $L^2(\mu)$.
\end{theorem*}
For $d=2$ the result was proved by Conze \& Lesigne~\cite{CL}. The 
particular case that the transformations $T_i$ are powers of the same 
transformation, for example $T_i= T^i$, was solved by Host \& 
Kra~\cite{HK1}.

Tao's proof does not  really belong to ergodic theory: he  uses only 
the pointwise ergodic theorem in order to translate the problem into a finitary 
question. Soon after, H. Townser~\cite{To} rewrote the proof using nonstandard 
analysis. More recently, T. Austin~\cite{A} gave another proof of the 
same 
result  by more conventional ergodic methods and the idea of the 
present work was inspired by  the reading of his paper.

Let us say a few words about the methods. All the papers dealing 
with a single transformation use an idea introduced by 
Furstenberg~\cite{F}:
the construction of a \emph{characteristic factor}. 
It is a factor (i.e. a quotient) of 
the system controlling the asymptotic behaviour of the multiple averages 
in a way that allow one to consider only functions defined on this 
factor. The next step is to prove that this factor has a nice 
\emph{structure}, and the convergence is much easier to prove in this 
case. 
In short, the convergence follows from the existence of a hidden 
structure of the system. The same structure can be used to study other 
problems of multiple convergence and of multiple recurrence, for 
example in 
\cite{HK2}, \cite{L}, \cite{BHK}, \cite{FHK}, \cite{FK1}, \dots \
A similar method was used by Conze \& Lesigne for two commuting 
transformations, but all attempts to solve the general case by using the 
machinery of characteristic factors  
 were unsuccessful\footnote{The case considered in~\cite{FK} is very 
particular.}.

T. Austin proceeds in the opposite 
direction, building an \emph{extension} of the original system with 
good properties; he calls it a \emph{pleasant system}. It happens 
that this 
extension is not very explicit (it is defined as an inverse limit) and 
that it gives little information  about the original system. Moreover,  its 
construction is directly related to the averages~\eqref{eq:av} and 
apparently is difficult to use  for related problems.

Although we give here a fourth proof of Tao's result, this is not the 
main concern of this paper. Our main goal is to provide some tools 
for the case of several commuting transformations,
similar to the tools 
successfully used in the case of a single transformation, 
with the idea that they will be 
useful in the solution of other problems. For this reason, we 
conclude this paper by adding in
Section~\ref{sec:permut} some properties that are we do not 
immediately need.

The price to pay for more generality is that some proofs in this paper are less 
elementary than in Austin's.

\subsection{}
Tao's method gives the convergence of the ordinary 
averages~\eqref{eq:av} only, while Austin's proof as well as ours 
generalizes to ``uniform averages'':
\begin{theorem}[T.~Austin]
\label{th:austin}
Let $(X,\mu,S_1,\dots,S_d)$ be a system where $S_1,\dots,S_d$ are 
commuting measure preserving transformations. \\ Then, for every 
$f_1,\dots,f_d\in L^\infty(\mu)$, the averages
\begin{equation}
\label{eq:av2}
\frac 1{|I_j|}\sum_{n\in I_j} S_1^nf_1\cdot\ldots\cdot S_d^nf_d
\end{equation}
converge in $L^2(\mu)$ for any sequence $(I_j\colon j\geq 1)$ of 
intervals in $\Z$ whose lengths $|I_j|$ tend to infinity.
\end{theorem}
In fact, Austin's result is slightly more general: instead of 
commuting transformations he considers commuting measure preserving 
$\Z^r$-actions on $X$;  the averages on intervals are replaced by 
averages on a F\o lner sequence in $\Z^r$. 
Up to minor changes (almost only in notation), the method presented here 
can be used  in this more general situation but for simplicity we 
restrict ourselves to the case stated in Theorem~\ref{th:austin}.

\subsection{Contents}
We first follow the same strategy as in the first sections 
of~\cite{HK1}:
Given a system $(X,\mu,T_1,\dots,T_d)$ where  the 
transformations commute, we build in Section~\ref{sec:construction} a measure
 $\mu^*$ on some Cartesian 
(finite) power $X^*$ of $X$ and use it to define a seminorm $\nnorm\cdot$ on 
$L^\infty(\mu)$ and we establish the properties of 
 that are used in the proof of 
Tao's Theorem. We show:

\begin{proposition}
\label{prop:caracteristic}
 Let $(X,\mu,S_1,\dots,S_d)$ be a system where $S_1,\dots,S_d$ are 
commuting measure preserving transformations. Define $T_1=S_1$ and 
$T_i=S_iS_1\inv$ for $2\leq i\leq d$ and let $\nnorm\cdot$ 
denote the seminorm on $L^\infty(\mu)$ associated to the system 
$(X,\mu,T_d,\dots,T_2,T_1)$.

 Then, for every 
$f_1,\dots,f_d\in L^\infty(\mu)$ with $\norm{f_i}_{L^\infty(\mu)}\leq 
1$ for $2\leq i\leq d$, we have
$$
 \limsup_{j\to+\infty}\Bigl\Vert 
\frac 1{|I_j|}\sum_{n\in I_j} S_1^nf_1\cdot\ldots\cdot 
S_d^nf_d\Bigr\Vert_{L^2(\mu)}\leq \nnorm{f_1}
$$
for every sequence  of intervals 
$(I_j\colon j\geq 1)$ in $\Z$ whose lengths tend to infinity.
\end{proposition}

Next, we remark that $X^*$ is naturally endowed with some commuting 
transformations $T_1^*,\dots,T_d^*$ and that $X^*$, endowed with 
$\mu^*$ and with these transformations, admits $X$ as a factor. 
Therefore, in order to prove the convergence of the 
averages~\eqref{eq:av}, we can substitute $X^*$ for $X$.

Properties of this  system are established in 
Section~\ref{sec:pretty}.
Substituting $(X^*,\mu^*,T_1^*,\dots,T_d^*)$ for 
$(X,\mu,T_1,$ $\dots,T_d)$, we define a seminorm $\nnorm\cdot^*$ on 
$L^\infty(\mu^*)$. The main result of this paper is:
\begin{theorem}
\label{th:pretty}
Let $\CW^*$ be the $\sigma$-algebra
$$
 \CW^*:=\bigvee_{i=1}^d \CI(T_i^*)
$$
of $(X^*,\mu^*)$, where $\CI(T_i^*)$ is the $\sigma$-algebra of sets 
invariant under $T_i^*$.

 If $F\in L^\infty(\mu^*)$ is such that 
$\E_{\mu^*}(F\mid\CW^*)=0$ then $\nnorm F^*=0$.
\end{theorem}
We call a system with this property a \emph{magic system}. 

Theorem~\ref{th:pretty} implies in particular that every system 
has a magic extension. 
This notion is similar to that of a pleasant system in~\cite{A} and is 
used in the same way. The 
differences are that $X^*$ is a relatively explicit system\footnote{It seems possible that the methods used 
here in the proof of Theorem~\ref{th:pretty} can be combined 
with the  constructions of~\cite{A}, removing the need for the inverse 
limit.} ($X^*$ is 
a finite cartesian power of $X$) and that its construction is related 
to the seminorm associated to the transformations and not only to the 
averages~\eqref{eq:av}. Therefore  it can be used to study any other 
question involving this seminorm.

 Tao's ergodic theorem follows 
easily from the preceding two results.

\begin{proof}[Proof of Theorem~\ref{th:austin}, assuming 
everything above] 
 By induction on $d$. For $d=1$, the statement is the mean 
ergodic theorem. We take $d>1$ and assume that the result is 
established for $d-1$ transformations. 

Let $T_1,\dots,T_d$ and $\nnorm\cdot$ be as in 
Proposition~\ref{prop:caracteristic}, 
$(X^*,\mu^*,T_1^*,\dots,T_d^*)$  as above and $\CW^*$  as in 
Theorem~\ref{th:pretty}. We define the transformations 
$S_1^*,\dots,S_d^*$ 
of $X^*$ by  $S_1^*=T_1^*$ and $S_i^*=T_i^*T_1^{*-1}$ for $2\leq 
i\leq d$. 

We have that  $(X,\mu,S_1,\dots,S_d)$ is 
a factor of $(X^*,\mu^*, S_1^*,\dots,S_d^*)$. Therefore, in order to 
prove the convergence of the averages~\eqref{eq:av2} in $L^2(\mu)$ for functions 
$f_1,\dots,f_d$ in $L^\infty(\mu)$, it suffices to show the 
convergence in $L^2(\mu^*)$ of the averages
\begin{equation}
\label{eq:avbis}
\frac 1{|I_j|}\sum_{n\in I_j} S_1^{*n}f_1^*\cdot\ldots\cdot S_d^{*n}f^*_d
\end{equation}
for functions $f_1^*,\dots,f_d^*$ in $L^\infty(\mu^*)$.

Consider first the case that 
\begin{equation}
\label{eq:decoposable}
 f_1^*=g_2\cdots g_d\text{ where } g_i\text{ is invariant under }
T_i^*\text{ for }2\leq i\leq d\ .
\end{equation}
As $T_i^*=S_i^*S_1^{*-1}$ for $2\leq i\leq d$ the 
averages~\eqref{eq:avbis} can be rewritten as 
$$
\frac 1{|I_j|}\sum_{n\in I_j} S_2^{*n}(g_2f_2^*)\cdot\ldots\cdot 
S_d^{*n}(g_df^*_d)
$$
and the convergence in $L^2(\mu^*)$ follows from the induction 
hypothesis.

Since the linear span of the functions of the 
form~\eqref{eq:decoposable} is dense in $L^\infty(\mu^*,\CW^*)$ for 
the norm of $L^1(\mu^*)$, we get by density that the 
averages~\eqref{eq:avbis} converge whenever the function $f_1^*$ is 
measurable with respect to $\CW^*$.

We are left with checking the case that $\E_{\mu^*}(f_1^*\mid\CW^*)=0$. 
 We have 
$\nnorm{f_1^*}^*=0$ by Theorem~\ref{th:pretty} and the 
averages~\eqref{eq:avbis} converge to $0$ in $L^2(\mu^*)$ by 
Proposition~\ref{prop:caracteristic}. 
\end{proof}

\section{The box measure and the box seminorm}
\label{sec:construction}
The objects defined in this section, as well as their properties, 
 are completely similar 
to those of Section~3 of~\cite{HK1}. Most of the proofs 
are exactly the same and we only sketch them.

\subsection{Notation} All functions are implicitly assumed to be 
measurable and real valued.

If $S$ is a measure preserving transformation of a probability space
$(Y,\nu)$, we write $\CI(S)$ for the algebra of  $S$-invariant sets. 
The \emph{conditionally independent square} of $\nu$ over $\CI(S)$ is the 
measure $\nu\times_{\CI(S)}\nu$ on $Y\times Y$ characterized by:\\
\emph{For all bounded measurable functions $f,f'$ on $X$,}
$$
 \int f(y)f'(y')\,d\nu\times_{\CI(S)}\nu\,(y,y')=\int 
\E_\nu\bigl(f\mid\CI(S)\bigr)\,\E_\nu\bigl(f'\mid\CI(S)\bigr)\,d\nu\ .
$$

We write $X^*=X^{2^d}$. We introduce some conventions for
notation of points in this space and more genrally in $X^{2^k}$ where 
$k\geq 1$ is an integer.

The points of $X^{2^k}$ are written
$$
  x=(x_\epsilon\colon \epsilon\in\{0,1\}^k)\ .
$$
Each $\epsilon\in\{0,1\}^k$ is written without commas and parentheses.
If $k\geq 2$ and $\eta\in\{0,1\}^{k-1}$, we write 
$\eta0=\eta_1\dots\eta_{k-1}0$ and $\eta 1=\eta_1\dots\eta_{k-1}1$.

Occasionally, it is  convenient to  also use another notation.
 We write $[k]=\{1,2,\dots,k\}$ and  make the natural identification 
between $\{0,1\}^k$ and the family of subsets of $[k]$. Therefore, 
for $\epsilon\in\{0,1\}^k$ and $1\leq i\leq k$, the assertion 
``$\epsilon_i=1$'' is equivalent to ``$i\in\epsilon$''. 
Therefore we  write  
$\emptyset=00\dots 0\in\{0,1\}^k$.

If $f_\epsilon$, $\epsilon\in\{0,1\}^k$, are functions on $X$, we 
define a function on $X^{2^k}$ by
$$
 \Bigl(\bigotimes_{\epsilon\in\{0,1\}^k}f_\epsilon\Bigr)(x)
:=\prod_{\epsilon\in\{0,1\}^k}f_\epsilon(x_\epsilon)\ .
$$

For $1\leq i\leq d$, $T_i^\Delta$ denotes the 
\emph{diagonal transformation}  
$T_i\times T_i\times\dots\times T_i$
of $X^{2^k}$:
$$
\text{ for every }
\epsilon\in\{0,1\}^d,\quad \bigl(T_i^\Delta x\bigr)_\epsilon=T_ix_\epsilon
$$
and the \emph{side transformations} $T_i^*$ of $X^*$ are given by
\begin{equation}
\label{eq:side}
\text{ for every }
\epsilon\in\{0,1\}^d,\quad
 (T^*_ix)_\epsilon=\begin{cases}
T_ix_\epsilon &\text{if }\epsilon_i=0\ ;\\
x_\epsilon &\text{if } \epsilon_i=1\ .
\end{cases}
\end{equation}

\subsection{The box measure}

We build a measure $\mu^*$ on $X^*$. First we define a measure 
$\mu_{T_1}$ on $X^2$ by 
$$\mu_{T_1}=\mu\times_{\CI(T_1)}\mu\ .$$
This means that for $f_0,f_1\in L^\infty(\mu)$ we have
\begin{equation}
\label{eq:mu1}
 \int f_0(x_0)f_1(x_1)\,d\mu_{T_1}(x)
=\int\E\bigl(f_0\mid\CI(T_1)\bigr)\cdot\E\bigl(f_1\mid\CI(T_1)\bigr)\,d\mu
\ .
\end{equation}
This measure is invariant under the transformations
$$
 T_i\times T_i\ (1\leq i\leq d)\text{ and }T_1\times \id\ .
$$
Next we define the measure $\mu_{T_1,T_2}$ on $X^4=X^2\times X^2$ by
$$
 \mu_{T_1,T_2}=\mu_{T_1}\times_{\CI(T_2\times T_2)}\mu_{T_1}\ .
$$
This means that for $f_{00},\dots,f_{11}\in L^\infty(\mu)$ we have
\begin{multline*}
\int \prod_{\epsilon\in\{0,1\}^2}f_\epsilon(x_\epsilon)\,d\mu_{T_1,T_2}(x)\\
=\int 
\E_{\mu_{T_1}}\bigl( f_{00}\otimes f_{10}\mid \CI(T_2\times T_2)\bigr)\cdot 
\E_{\mu_{T_1}}\bigl( f_{01}\otimes f_{11}\mid \CI(T_2\times 
T_2)\bigr)\,d\mu_{T_1}
\ .
\end{multline*}
For $1\leq i\leq d$, this measure is invariant under the ``diagonal 
transformations'' 
$ T_i\times T_i\times T_i\times T_i$ of $X^4$; it is also invariant 
under the ``side transformations''  
$T_1\times\id\times T_1\times\id$ and $ T_2\times 
T_2\times\id\times\id$.

In the same way, for $k<d$ we obtain a measure $\mu_{T_1,\dots,T_k}$ on $X^{2^k}$, 
invariant under  all ``diagonal transformations'' 
$T_i\times T_i\times\dots\times T_i$ ($1\leq i\leq d$) and under the 
``side 
transformations'' associated to $T_1,\dots,T_k$ as in~\eqref{eq:side}, 
but with $k$ substituted for $d$. We define:
$$
 \mu_{T_1,\dots,T_{k+1}}=\mu_{T_1,\dots,T_k}
\times_{\CI(T_{k+1}\times T_{k+1}\times\dots\times 
T_{k+1})}\mu_{T_1,\dots,T_k}\ .
$$

After $d$ steps we obtain a measure
 $\mu^*:=\mu_{T_1,\dots,T_d}$ on $X^*=X^{2^d}$.
If $f_\epsilon$, $\epsilon\in\{0,1\}^d$, belong to 
$L^\infty(\mu)$, we have
\begin{multline}
\label{eq:iterate}
\int 
\bigotimes_{\epsilon\in\{0,1\}^d}f_\epsilon\,d\mu^*(x)\\
=\int\E_{\mu_{T_1,\dots,T_{d-1}}}\Bigl(\bigotimes_{\eta\in\{0,1\}^{d-1}}
f_{\eta 0}\big\vert\CI(T_d\times\dots\times T_d)\Bigr)\quad\strut\\
\strut\hskip 3cm \cdot 
\E_{\mu_{T_1,\dots,T_{d-1}}}\Bigl(\bigotimes_{\eta\in\{0,1\}^{d-1}}f_{\eta 
1}\big\vert \CI(T_d\times\dots\times 
T_d)\Bigr)\,d\mu_{T_1,\dots,T_{d-1}}
\end{multline}
and thus 
\begin{equation}
\label{eq:limd}
\int 
\bigotimes_{\epsilon\in\{0,1\}^d}f_\epsilon\,d\mu^*(x)
=\lim_{N\to+\infty}\frac 1{N}\sum_{n=0}^{N-1}\int
\bigotimes_{\eta\in\{0,1\}^{d-1}}\bigl(T_d^{n}f_{\eta 0}\cdot 
f_{\eta_1}\bigr)\,d\mu_{T_1,\dots,T_{d-1}}\ .
\end{equation} 
Moreover, the same convergence holds if the intervals $[0,N)$ are 
replaced by  any sequence of intervals of lengths tending to infinity.
Starting from~\eqref{eq:limd} and proceding by downwards induction we 
get:
\begin{lemma}
\label{lem:limit}
 If $f_\epsilon$, $\epsilon\in\{0,1\}^d$ belong to 
$L^\infty(\mu)$, we have
\begin{multline}
\label{eq:limit}
 \int \prod_{\epsilon\in\{0,1\}^d} f_\epsilon(x_\epsilon)\,d\mu^*(x)\\
=
\lim_{N_d\to+\infty}\frac 1{N_d}\sum_{n_d=0}^{N_d-1}\dots
\lim_{N_2\to+\infty}\frac 1{N_2}\sum_{n_2=0}^ {N_2-1}
\lim_{N_1\to+\infty}\frac 1{N_1}\sum_{n_1=0}^ {N_1-1}\\
\int 
\prod_{\epsilon\in\{0,1\}^d}T_1^{(1-\epsilon_1)n_1}\cdots
T_d^{(1-\epsilon_d)n_d}f_\epsilon\,d\mu\ .
\end{multline}

Moreover, relation~\eqref{eq:limit} holds for averages on any other sequence of 
intervals whose length tends to infinity, 
for example for the symmetric averages on $[-N_i,N_i]$.
\end{lemma}

The  measure $\mu^*$ is
invariant under the diagonal transformations $T_i^\Delta$ and the 
side transformations  $T^*_i$, $1\leq i\leq d$.
This measure is called the \emph{box measure associated to the 
transformations $T_1,\dots,T_d$}.

In some cases we write $\mu_{T_1,\dots,T_d}$ instead of $\mu^*$ to avoid any possible 
ambiguity.

We notice that all the marginals of $\mu^*$ are equal to $\mu$ and 
that the projection $\pi_\emptyset\colon X^{2^k}\to X$ given by 
$\pi_\emptyset(x)=x_\emptyset$ 
is a factor map from $(X^*,\mu^*, T_1^*,\dots,T_d^*)$ to 
$(X,\mu,T_1,\dots,T_d)$.

For $1\leq i\leq d$  the coordinate indexed by any
$\epsilon\in\{0,1\}^d$ plays the same role in the 
construction of $\mu^*$ as the coordinate indexed by $\epsilon'$ 
obtained in substituting $1-\epsilon_i$ for $\epsilon_i$. This shows 
that the measure $\mu^*$ is invariant under the symmetry of $X^*$ 
associated in the obvious way to this map.

\subsection{The box seminorm}
By~\eqref{eq:iterate}, 
 for every $f\in L^\infty(\mu)$ we have
$$
 \int\prod_{\epsilon\in\{0,1\} ^d} 
f(x_\epsilon)\,d\mu^*(x)
\geq 0
$$
and we can define:
\begin{definition}
For $f\in L^\infty(\mu)$,
\begin{equation}
\label{eq:defnorm}
 \nnorm f:=\Bigl(
 \int\prod_{\epsilon\in\{0,1\}^d} f(x_\epsilon)\,d\mu^*(x)\Bigr)^{1/2^d}\ .
\end{equation}
When needed we write $\nnorm f_{T_1,\dots,T_d}$ instead of $\nnorm f$.
\end{definition}
From~\eqref{eq:limd} we get:

\emph{For every $f\in L^\infty(\mu)$ we have}
\begin{equation}
\label{eq:recurnorm}
\nnorm f_{T_1,\dots,T_d}^{2^d}=\lim_{N_d\to+\infty}
\sum_{n_d=0}^{N_d-1}
\nnorm{T_d^{n_d}f\cdot  f}_{T_1,\dots,T_{d-1}}^{2^{d-1}}\ .
\end{equation}

\begin{remark}
As in~\cite{HK1}, a similar formula can be derived for complex valued 
functions. We do not give it here.
\end{remark}

\begin{proposition}[and definition]
\label{prop:seminorm}\strut 
\begin{enumerate}
\item
For $f_\epsilon\in L^\infty(\mu)$, $\epsilon\in\{0,1\}^d$, we have
\begin{equation}
\label{eq:CSG}
\Bigl|\int \bigotimes_{\epsilon\in\{0,1\}^d} 
f_\epsilon\,d\mu^*\Bigr|\leq
\prod_{\epsilon\in\{0,1\}^d} \nnorm{f_\epsilon}\ .
\end{equation}
\item $\nnorm\cdot$ is a seminorm on $L^\infty(\mu)$.\\
We call it \emph{the box seminorm} associated to $T_1,\dots,T_d$.
\end{enumerate}
\end{proposition}

The bound~\eqref{eq:CSG} is similar to the Cauchy-Schwarz-Gowers 
Inequality.
\begin{proof}
 The first part of the Proposition is proved by induction on 
$d$. For $d=1$, the result follows imediately from the 
definition~\eqref{eq:mu1} of $\mu_{T_1}$ and 
the Cauchy-Schwarz Inequality. We assume now that $d\geq 2$ and that 
the result is true for $d-1$ transformations.

For $\epsilon\in\{0,1\}^d$ we define two functions $f'_\epsilon$ and 
$f''_\epsilon$ on $X$ by
$$
\text{for all }\eta\in\{0,1\}^{d-1},\quad
 f'_{\eta0}=f'_{\eta1}=f_{\eta0}\text{ and }
f''_{\eta0}=f''_{\eta1}=f_{\eta1}\ .
$$
Let $I$ be the left hand side of~\eqref{eq:CSG} and let $I'$ and $I''$ 
be respectively the similar expressions obtained by substituting the 
functions $f'_\epsilon$, respectively $f''_\epsilon$,  for the 
functions $f_\epsilon$.
By~\eqref{eq:iterate} and the Cauchy-Schwarz Inequality, $I^2\leq 
I'I''$.

By~\eqref{eq:limd}, the induction hypothesis, H\"older Inequality 
and~\eqref{eq:recurnorm},
\begin{multline*}
 I'=\Bigl\vert\lim_{N_d\to+\infty}\frac 1{N_d}\sum_{n_d=0}^{N_d-1}\int
\bigotimes_{\eta\in\{0,1\}^{d-1}}\bigl(T_d^nf_{\eta 0}\cdot 
f_{\eta_0}\bigr)\,d\mu_{T_1,\dots,T_{d-1}}\bigr\vert\\
\leq\limsup_{N_d\to+\infty}\frac 1{N_d}\sum_{n_d=0}^{N_d-1}
\prod_{\eta\in\{0,1\}^{d-1}} 
\nnorm{T_d^{n_d}f_{\eta 0}\cdot f_{\eta 0}}_{T_1,\dots,T_{d-1}}\\
\leq\prod_{\eta\in\{0,1\}^{d-1}}
\nnorm{f_{\eta 0}}_{T_1,\dots,T_d}^2
\ .
\end{multline*}
A similar bound holds for $I''$ and the result follows.

The second part of the proposition is obtained by using the same proof as 
for Lemma~3.9 in~\cite{HK1}.

\end{proof}

\subsection{Proof of Proposition~\ref{prop:caracteristic}}
The proof is the same as that of results for a single 
transformation, for example of Theorem 12.1 of~\cite{HK1}.
The proof proceeds by induction on $d$. For $d=1$ the seminorm is the 
absolute value of the integral and there is nothing to prove. We set 
$d>1$ and assume that the result is true for $d-1$ transformations.

Let $f_1,\dots,f_d$ and $S_1,\dots,S_d$ be as in the proposition.
We recall that $T_1=S_1$ and that $T_i=S_iS_1\inv$ for $2\leq i\leq d$.
By the van der Corput Lemma and Cauchy-Schwarz Inequality, the $\limsup$ in the 
proposition is bounded by 
$$
 \limsup_{H\to+\infty}\frac 1H\sum_{h=0}^{H-1} \limsup_{j\to+\infty}
\Bigl\Vert\frac1{|I_j|} \sum_{n\in I_j}\prod_{\substack{1\leq i\leq d\\ 
i\neq 2}}
(S_iS_2\inv)^n\bigl( f_i\cdot 
S_i^hf_i\bigr)\Bigr\Vert_{L^2(\mu)}\ .
$$

By the induction hypothesis, 
this  $\limsup$ is bounded by 
$$
 \limsup_{H\to+\infty}\frac 1H\sum_{h=0}^{H-1}\nnorm{ f_1\cdot 
S_1^hf_1}^\sharp\ ,
$$
where 
$\nnorm\cdot^\sharp$ is the seminorm associated to the 
transformations $(S_dS_2\inv)(S_1S_2\inv)\inv=T_d$,\dots, $
(S_3S_2\inv)(S_1S_2\inv)\inv=T_3$ and $S_1S_2\inv=T_2\inv$.
By construction, this seminorm remains unchanged if $T_2$ is 
substituted for $T_2\inv$ and thus is equal to the seminorm 
$\nnorm\cdot_{T_d,\dots,T_3,T_2}$.

 By Lemma~\ref{lem:limit} and Corollary~\ref{cor:permut},
$$
 \frac 1H\sum_{h=0}^{H-1}\nnorm{f_1\cdot 
S_1^hf_1}^{
2^{d-1}}_{T_d,\dots,T_2}\to\nnorm{f_1}^{2^d}_{T_d,\dots,T_2,S_1}\text{ as }H\to+\infty
$$
and we are done since $S_1=T_1$.\qed

\subsection{A uniformity result}
\label{subsec:uniform}
The next Lemma has no analogue in~\cite{HK1}.

\begin{lemma}
\label{lem:uniform}
Let $f_\emptyset\in L^\infty(\mu)$. Then for 
every $\delta>0$ there exists $N_0=N_0(\delta)$ such that:\\
For all $f_\epsilon\in L^\infty(\mu)$,
 $\emptyset\neq \epsilon\in\{0,1\}^d$ 
with $\norm{f_\epsilon}_{L^\infty(\mu)}\leq 1$, 
for all intervals $I_1,\dots,I_d$ of $\Z$ of length $\geq N_0$,
$$
\Bigl|\,\frac 1{|I_1|\cdot\dots\cdot|I_d|}\sum_{\substack{n_1\in 
I_1\\\dots\\n_d\in 
I_d}}
\int\prod_{\epsilon\in\{0,1\}^d}T_1^{(1-\epsilon_1)n_1}\dots 
T_d^{(1-\epsilon_d)n_d}f_\epsilon\,d\mu\,\Bigr|<\nnorm{f_\emptyset}+\delta\ .
$$
\end{lemma}

\begin{proof} 
We can assume that $\norm{f_\emptyset}_{L^\infty(\mu)}\leq 1$. 

Let $J$ 
be the average in the statement and
let $H_1,\dots,H_d$ be integers with $1\leq H_i\leq|I_i|$ for all $i$.

Each $\epsilon\in\{0,1\}^d$ is written either $\epsilon=\eta 0$ with 
$\eta\in\{0,1\}^{d-1}$ or 
$\epsilon=\eta 1$, depending on the value  of $\epsilon_d$.
We split the product in the integral in two parts: 
\begin{enumerate}
\item
The product of the terms indexed by $\eta 0$ for some 
$\eta\in\{0,1\}^{d-1}$. This product can be written as 
 $T_d^{n_d}F_{n_1,\dots,n_{d-1}}$.
\item
The product $F'_{n_1,\dots,n_{d-1}}$ of the terms indexed by 
$\eta 1$ for some 
$\eta\in\{0,1\}^{d-1}$. 
\end{enumerate}
We thus have that $J$ is equal to
$$
\frac 1{|I_1|\cdot\dots\cdot|I_{d-1}|}
\sum_{\substack{n_1\in I_1 \\\dots\\ n_{d-1}\in I_{d-1}}}
\int \frac 1{|I_d|}\sum_{n_d\in I_d}
T_d^{n_d}F_{n_1,\dots,n_{d-1}}\cdot F'_{n_1,\dots,n_{d-1}}\, d\mu
$$
and as $|F'_{n_1,\dots,n_{d-1}}|\leq 1$ we have
$$
 |J|^2\leq \frac 1{|I_1|\cdot\dots\cdot|I_{d-1}|}
\sum_{\substack{n_1\in I_1 \\\dots\\ n_{d-1}\in I_{d-1}}}
\Bigl\Vert 
 \frac 1{|I_d|}\sum_{n_d\in I_d}
T_d^{n_d}F_{n_1,\dots,n_{d-1}}
\Bigr\Vert_{L^2(\mu)}^2\ .
$$
By the finite van der Corput Lemma, the square of 
the norm in this formula  is 
bounded by the absolute value of
$$
\frac{4H_d}{|I_d|}+
\sum_{h_d=-H_d}^{H_d}\frac{H_d-|h_d|}{H_d^2}
\int T_d^{h_d} F_{n_1,\dots,n_{d-1}}\cdot F_{n_1,\dots,n_{d-1}}\,
d\mu\ .
$$
Replacing $F$ by its value, we get that $|J|^2$ is bounded by the 
absolute value of
\begin{multline*}
\frac{4H_d}{|I_d|}+
\frac 1{|I_1|\cdot\dots\cdot|I_{d-1}|}
\sum_{\substack{n_1\in I_1 \\\dots\\ n_{d-1}\in I_{d-1}}}
 \sum_{h_d=-H_d}^{H_d}\frac{H_d-|h_d|}{H_d^2}\\
\int \prod_{\epsilon\in\{0,1\}^d} 
T_1^{(1-\epsilon_1)n_1}\dots 
T_{d-1}^{(1-\epsilon_{d-1})n_{d-1}}
T_d^{(1-\epsilon_d)h_d}g_\epsilon\, d\mu
\end{multline*}
where the functions $g_\epsilon$ are given by 
$g_{\eta 0}=g_{\eta 1}=f_{\epsilon}$
 for $\eta\in\{0,1\}^{d-1}$.

We iterate the same computation, using successively  
$\epsilon_{d-1}\dots,\epsilon_2,\epsilon_1$ instead of 
$\epsilon_d$. We get that 
\begin{multline*}
 |J|^{2^d}\leq 
C\bigl(\frac{H_1}{|I_1|}+\dots+\frac{H_d}{|I_d|}\bigr)\\
+
\Bigl|\sum_{\substack{-H_1\leq h_1\leq H_1\\ \dots \\
-H_d\leq h_d\leq H_d}}
\prod_{i=1}^d\frac{H_i-|h_i|} {H_i^2}
\int\prod_{\epsilon\in\{0,1\}^d}
T_1^{(1-\epsilon_1)h_1}\dots T_d^{(1-\epsilon_d)h_d}
f_\emptyset\,d\mu\Bigr|
\end{multline*}
for some absolute constant $C$.

The iterated limit of the last average when $H_1\to+\infty$,\dots, 
$H_d\to+\infty$ is equal to $\nnorm{f_\emptyset}^{2^d}$ 
by Lemma~\ref{lem:limit}.
Therefore there exist $H_1,\dots,H_d$ such that this average has an 
absolute value less that $(\nnorm {f_\emptyset}+\delta/2)^{2^d}$. The result follows.
\end{proof}

\begin{remark}
\label{rem:permut}
It is easy to 
to check that the 
role played by $f_\emptyset$ in Lemma~\ref{lem:uniform} can be played 
by $f_\eta$ for any $\eta\in\{0,1\}^d$ and this implies a weak 
version of the bound~\eqref{eq:CSG} in 
Proposition~\ref{prop:seminorm}: the integral in the left hand member 
is equal to zero whenever at least one of the functions $f_\epsilon$ 
has  zero seminorm. In fact, this weak version would suffice 
or our purpose.
\end{remark}

\subsection{A characteristic $\sigma$-algebra on $X$}
The definitions and  results of this section are  completely similar to 
those of Section 4.2 of~\cite{HK1}.

Let us identify $X^*=X^{2^d}$ with $X^{2^{d-1}}\times X^{2^{d-1}}$; each 
point $x\in X^*$ is written $x=(x',x'')$, where $x',x''\in 
X^{2^{d-1}}$ are given by:
$$
x'=(x_{\eta 0}\colon \eta\in\{0,1\}^{d-1})\text{ and }
x''=(x_{\eta 1}\colon \eta\in\{0,1\}^{d-1})\ .
$$
 By construction, the 
images of $\mu^*$ under the projections $x\mapsto x'$ and $x\mapsto 
x''$ are equal to the measure $\mu_{d-1}$ 
associated to the transformations $T_1,\dots,T_{d-1}$. 
We remark also that 
\begin{equation}
\label{eq:invard}
 T_d^\Delta T_d^{*-1}=\id\times T_d^\circ\text{ where }T_d^\circ= \underbrace{T_d\times\dots\times T_d}
_{2^{d-1}\text{ times}}
\end{equation}

From the  inductive definition of the measure $\mu^*$, we deduce:

\begin{lemma}
\label{lem:invarTd}
Let $F\in L^\infty(\mu^*)$ be a function invariant under the 
transformation $T_d^\Delta T_d^{*-1}$. Then there exists a function 
$G$ on $X^{2^{d-1}}$, belonging to $L^\infty(\mu_{d-1})$, such 
that 
$$
F(x)=G(x')\quad \text{for }\mu^*\text{-almost every }x=(x',x'')\in X^*\ .
$$
\end{lemma}

By induction on $d$, we get:
\begin{corollary}\label{cor:invarTd}
Let $F\in L^\infty(\mu^*)$ be a function invariant under the 
transformations $T_i^\Delta T_i^{*-1}$ for $i=1,\dots,d$. Then there 
exists a function $f\in L^\infty(\mu)$ such that 
$F(x)=f(x_\emptyset)$ for $\mu^*$-almost every $x\in X^*$.
\end{corollary}

We write $X^\sharp=X^{2^d-1}$ and identify $X^*$ with $X\times 
X^\sharp$ by isolating the coordinate $\emptyset$ of each point: every 
point $x\in X^*$ is written
$$
 x=(x_\emptyset,x^\sharp)\text{ where }x^\sharp=
\bigl(x_\epsilon\colon \epsilon\in\{0,1\}^d,\ \epsilon\neq \emptyset\bigr)
\in X^\sharp\ .
$$
We write $\mu^\sharp$ for the image of $\mu^*$ in $X^\sharp$ under 
the  projection $x\mapsto x^\sharp$.

For $1\leq i\leq d$, the measure preserving  transformation 
$T_i^\Delta T_i^{*-1}$  of $(X^*,\mu^*)$ leaves 
the coordinate $x_\emptyset$ of each point $x$ invariant, and thus we 
can write this transformation as   
$$
T_i^\Delta T_i^{*-1}=\id_X\times  T_i^\sharp
$$ 
where $T_i^\sharp$ is the measure preserving 
transformation of $(X^\sharp,\mu^\sharp)$ given by
$$
\text{for }\emptyset\neq\epsilon\in\{0,1\}^d,\quad
 (T_i^\sharp x)_\epsilon=
\begin{cases}
T_ix_\epsilon &\text{ if }\epsilon_i=1\ ;\\
x_\epsilon &\text{ if }\epsilon_i=0 \ .
\end{cases}
$$

From Corollary~\ref{cor:invarTd} we immediately deduce:
\begin{corollary}
\label{cor:J}
Let $\CJ^\sharp$ be the $\sigma$-algebra of invariant sets of 
$(X^\sharp,\mu^\sharp,T_1^\sharp,\dots,$ $T_d^\sharp)$.

Then for every $A\in\CJ^\sharp$ there exists a subset $B$ of $X$ with
\begin{equation}
\label{eq:xxsharp}
\one_B(x_\emptyset)=\one_A(x^\sharp)\text{ for $\mu^*$-almost every 
}x=(x_\emptyset,x^\sharp)\in X^*\ .
\end{equation}
\end{corollary}
We remark that conversely, if $A\subset X^\sharp$ and $B\subset X$ 
satisfy~\eqref{eq:xxsharp}, then $A$ is invariant under $T_i^\sharp$ 
for every $i$.

\begin{lemma}[\cite{HK1}, Lemma~4.3]
\label{lem:Z}
Let $\CZ$ be the $\sigma$-algebra on $X$ consisting in sets $B$ such 
that there exists a subset $A$ of $X^\sharp$ satisfying the 
relation~\eqref{eq:xxsharp} of Corollary~\ref{cor:J}.

Then, for every $f\in L^\infty(\mu)$ we have\\
\centerline{ $\nnorm f=0$ if 
and only if $\E_\mu(f\mid\CZ)=0$.}
\end{lemma}

\begin{proof}
Assume first that $\E_\mu(f\mid\CZ)=0$. Let $F$ be the function on 
$X^\sharp$ given by
$$
 F(x^\sharp)=\prod_{\emptyset\neq \epsilon\in\{0,1\}^d} 
f(x_\epsilon)\ .
$$
Let $\CJ^\sharp$ be defined as in Corollary~\ref{cor:J}.
The function $x\mapsto 
\E_{\mu^\sharp}\bigl(F\mid\CJ^\sharp)\bigr)(x^\sharp)$  on $X^*$ 
is invariant under all 
transformations $T_i^\Delta T_i^{*-1}$ and thus by 
Corollary~\ref{cor:invarTd} there exists a function $g$ on $X$ with 
$$
 g(x_\emptyset)=\E_{\mu^\sharp}(F\mid\CJ^\sharp)(x^\sharp)
\text{ for }\mu^*\text{-almost every }
x=(x_\emptyset,x^\sharp)\ .
$$
 As $\mu^*$ is invariant under 
$\id_X\times  T_i^\sharp$ for every $i$, by definition of the 
seminorm we have
\begin{multline*}
\nnorm f^{2^d}=\int_{X^*} f(x_\emptyset) 
F(x^\sharp)\,d\mu^*(x_\emptyset,x^\sharp)\\
=\int_{X^*}  f(x_\emptyset)
\E_{\mu^\sharp}(F\mid\CJ^\sharp)(x^\sharp)\,d\mu^*(x)
=\int_X f(x_\emptyset) g(x_\emptyset)\,d\mu(x_\emptyset)=0
\end{multline*}
because $g$ is measurable with respect to $\CZ$ by definition.

We assume now that $\nnorm f=0$. Let $g\in L^\infty(\mu)$ be 
measurable with respect to $\CZ$. By definition, there exists a 
function $G\in L^\infty(\mu^\sharp)$ with $g(x_\emptyset)=G(x^\sharp)$, 
$\mu^*$-almost everywhere. We have
$$
 \int_X f(x)g(x)\,d\mu(x)=\int_{X^*} 
f(x_\emptyset)G(x^\sharp)\,d\mu^*(x_\emptyset,x^\sharp)
$$
and it follows from the bound~\eqref{eq:CSG} of Proposition~\ref{prop:seminorm}
that this integral is equal to zero.
\end{proof}

In the case of  single transformation, the $\sigma$-algebra $\CZ$ 
is the $\sigma$-algebra $\CZ_{d-1}$ of~\cite{HK1}, where it is shown 
that the corresponding factor $Z_{d-1}$ has the structure of an 
\emph{inverse limit of $(d-1)$-step nilsystems}. But in the present 
case of several transformations $\CZ$  apparently only has a weaker 
structure and
 we stop following~\cite{HK1} at this 
point.
\section{Proof of Theorem~\ref{th:pretty}}
\label{sec:pretty}

\subsection{The system $(X^*,\mu^*,T_1^*,\dots,T_d^*)$ }

Let $\CX^\sharp$ be the $\sigma$-algebra on $X^*$ corresponding to the 
factor $X^\sharp$ of $X^*$: $\CX^\sharp$ is spanned by the 
projections $x\mapsto x_\epsilon\colon X^*\to X$ for 
$\epsilon\in\{0,1\}^d$, $\epsilon\neq \emptyset$.

\begin{lemma}
\label{lem:span0}
The subspace of $L^2(\mu^*)$ consisting in functions with zero 
conditional expectation on $\CX^\sharp$ is the space spanned by 
functions of the 
form 
$$
F(x)=\prod_{\epsilon\in\{0,1\}^d}
f_\epsilon(x_\epsilon) \text{ where }|f_\epsilon|\leq 1\text{ 
for all }\epsilon\text{ and }\E_\mu(f_\emptyset\mid\CZ)=0\ .
$$
\end{lemma}

\begin{proof}
Let $\CL$ be the  closed subspace of $L^2(\mu^*)$ spanned by functions of the 
type given in the statement and let $\CL'$ be the  closed subspace of 
$L^2(\mu^*)$ spanned by functions of the form 
$$
 F'(x)=\prod_{\epsilon\in\{0,1\}^d}
f'_\epsilon(x_\epsilon) \text{ where }|f'_\epsilon|\leq 1\text{ 
for all }\epsilon\text{ and }f'_{\emptyset}\text{ is $\CZ$-measurable.}
$$

The sum of these spaces is clearly dense in 
$L^2(X^*,\mu^*)$. We claim that they are orthogonal.

Let $f_\epsilon$ and $f'_\epsilon$, $\epsilon\in\{0,1\}^d$, be as 
above.
For every $i$, the function $x\mapsto 
f_\emptyset(x_\emptyset)f'_\emptyset(x_\emptyset)$ is invariant under 
$\id_X\times T_i^\sharp$ and thus
$$
 \int f_\emptyset(x_\emptyset)f'_\emptyset(x_\emptyset)
\prod_{\emptyset\neq\epsilon\in\{0,1\}^d}f(x_\epsilon)f'(x_\epsilon)
\, d\mu^*
=\int f_\emptyset(x_\emptyset)f'_\emptyset(x_\emptyset) G(x^\sharp)\,d\mu^*
$$ where
$$
 G =\E_{\mu^\sharp}\Bigl(
\bigotimes_{\emptyset\neq\epsilon\in\{0,1\}^d}f_\epsilon 
f'_\epsilon\big\vert\CJ^\sharp\Bigr)\ .
$$
By Corollary~\ref{cor:J} there exists a function $g\in L^\infty(\mu)$, 
measurable with respect to $\CZ$, with $g(x_\emptyset)=G(x^\sharp)$ for 
$\mu^*$-almost every $x=(x_\emptyset,x^\sharp)$ and the integral above 
is equal to
$$
 \int f_\emptyset(x_\emptyset)f'_\emptyset(x_\emptyset)  
g(x_\emptyset)\,d\mu(x_\emptyset)\ .
$$
This is equal to zero because $\E_\mu(f_\emptyset\mid\CZ)=0$ and the 
function $f'_\emptyset g$ is measurable with respect to $\CZ$. Our claim 
is proved.
Therefore $\CL$ is the orthogonal space to $\CL'$.

On the other hand, $\CL'$ clearly contains 
$L^2(X^*,\CX^\sharp,\mu^*)$ and by the definition of $\CZ$ in 
Lemma~\ref{lem:Z} we have the opposite inclusion and so
these spaces are equal. Therefore, $\CL$ is the orthogonal space to
$L^2(X^*,\CX^\sharp,\mu^*)$, and this is the announced result.
\end{proof}

\subsection{Iterating the construction}\strut\\
We now define a new system $(X^{**},\mu^{**},T_1^{**},\dots, T_d^{**})$ where 
$X^{**}:=(X^*)^*=(X^{2^d})^{2^d}$. 
It is built from the system $(X^*,\mu^*,T_1^*,\dots, T_d^*)$
  in the same way that  
$(X^*,\mu^*,T_1^*,\dots,T_d^*)$ was built from $(X,\mu,T_1,$ $\dots,T_d)$.
The points of $X^{**}$ are written 
$$
x=\bigl( x_{\epsilon\eta}\colon\epsilon,\eta\in\{0,1\}^d\bigr)\ ,
$$
with the $2^d$ natural projections $\pi_\eta^*\colon X^{**}\to X^*$ being given by the maps 
$$
\bigl( \pi_\eta^*(x)\bigr)_\epsilon=x_{\epsilon\eta}\ .
$$
The seminorm $\nnorm\cdot^*$ on 
$L^\infty(\mu^*)$ is defined from the measure $\mu^{**}$ in the same 
way as the seminorm $\nnorm\cdot$ on $L^\infty(\mu)$ was defined from the measure 
$\mu ^*$. 
\begin{lemma}\label{lem:normstar}
Let
 $$F(x)=\prod_{\epsilon\in\{0,1\}^d}f_\epsilon(x_\epsilon)\text{ 
where }f_\epsilon\in L^\infty(\mu)\text{ for all }\epsilon\text{ and 
}\nnorm {f_\emptyset} =0\ .
$$
Then $\nnorm F^*=0$.
\end{lemma}

\begin{proof}
We can assume that $|f_\epsilon|\leq 1$ for $\epsilon\neq\emptyset$.
By Lemma~\ref{lem:limit} applied to the measure $\mu^{**}$, 
$\nnorm F^{*2^d}$ is equal to the iterated
limit when $P_1,\dots,P_d\to+\infty$ of the averages for 
$p_1\in[0,P_1)$, \dots,$p_d\in[0,P_d)$ of 
$$
 I(p_1,\dots,p_d):=\int \prod_{\eta\in\{0,1\}^d} T_1^{*\; (1-\eta_1)p_1}\dots 
T_d^{*\;(1-\eta_d)p_d}\Bigl(\bigotimes_{\epsilon \in\{0,1\}^d}
f_\epsilon\Bigr)\,d\mu^*\ .
$$
By  definition of the transformations $T_i^*$, this is equal 
to
$$
 \int \bigotimes_{\epsilon \in\{0,1\}^d}\Bigl(\prod_{\eta\in\{0,1\}^d}
T_1^{(1-\eta_1)(1-\epsilon_1)p_1}\dots 
T_d^{(1-\eta_d)(1-\epsilon_d)p_d}f_{\epsilon}\Bigr)\,d\mu^*\ .
$$
By Lemma~\ref{lem:limit} again, but now applied to the measure $\mu ^*$,
 $\nnorm F^{*2^d}$ is equal to the iterated
limit when $N_1,\dots,N_d,P_1,\dots,P_d\to+\infty$ of the averages for 
$n_1\in[0,N_1)$,\dots, $n_d\in[0,N_d)$, $p_1\in[0,P_1)$, \dots,$p_d\in[0,P_d)$ of 
\begin{multline*}
J(n_1,\dots,n_d,p_1,\dots,p_d):=\\
\int 
 \prod_{\epsilon,\eta\in\{0,1\}^d}
T_1^{(1-\epsilon_1)(1-\eta_1)p_1+(1-\epsilon_1)n_1}\dots
T_d^{(1-\epsilon_d)(1-\eta_d)p_d+(1-\epsilon_d)n_d)}
f_{\epsilon}\,d\mu\ .
\end{multline*}

At this point, it is more convenient to identify $\{0,1\}^d$ 
with the family of subsets 
of $[d]$. Let $\theta\subset[d]$. In the product in $\epsilon,\eta$ of the 
last formula,  we gather all the terms 
 with $\epsilon\cup\eta=\theta$.  For $1\leq i\leq d$ we have 
$(1-\epsilon_i)(1-\eta_i)p_i+(1-\epsilon_i)n_i=
(1-\theta_i)(p_i+n_i)+\eta_in_i$.
We get that
\begin{multline*}
J(n_1,\dots,n_d,p_1,\dots,p_d)\\
=\int 
 \prod_{\theta\subset[d]}
T_1^{(1-\theta_1)(p_1+n_1)}\dots 
T_d^{(1-\theta_d)(p_d+n_d)}g_\theta^{(n_1,\dots,n_d)}\,d\mu
\end{multline*}
where
$$
 g_\theta^{(n_1,\dots,n_d)}=\prod_{\eta\subset\theta}
T_1^{\eta_1n_1}\dots T_d^{\eta_dn_d} 
\prod_{\epsilon\colon\epsilon\cup\eta=\theta} f_\epsilon\ .
$$
We consider $P_1,\dots,P_d$ as fixed. We have:
\begin{multline*}
K(n_1,\dots,n_d)
:= \frac 1{P_1\dots 
P_d}\sum_{p_1=0}^{P_1-1}\dots\sum_{p_d=0}^{P_d-1}
J(n_1,\dots,n_d,p_1,\dots,p_d) 
\\
= \frac 1{P_1\dots 
P_d}\sum_{p_1=n_1}^{n_1+P_1-1}\dots\sum_{p_d=n_d}^{n_d+P_d-1}
\\
\int   \prod_{\theta\subset[d]}
T_1^{(1-\theta_1)p_1}\dots 
T_d^{(1-\theta_d)p_d}g_\theta^{(n_1,\dots,n_d)}\,d\mu\ .
\end{multline*}
We remark that for every $n_1,\dots,n_d$ we have
$$
 |g_\theta^{(n_1,\dots,n_d)} |\leq 1\text{ for every $\theta$ and }
g_\emptyset^{(n_1,\dots,n_d)}=f_\emptyset\ .
$$
Therefore, by Lemma~\ref{lem:uniform}, 
for every $\delta>0$ there exists $P$ such that
$$
 |K(n_1,\dots,n_d)|<\delta\text{ for all }n_1,\dots,n_d\text{ 
whenever } P_1,\dots,P_d>P
$$
and the announced conclusion follows.
\end{proof}

\subsection{End of the proof}
We recall that $\CW^*$ is the $\sigma$-algebra 
$$
 \CW^*=\bigvee_{i=1}^d\CI(T_i^*)
$$
 on $(X^*,\mu^*)$. We show that if a 
function $F\in L^\infty(\mu^*)$ satisfies $\E_{\mu^*}(F\mid\CW^*)=0$ then 
$\nnorm F^*=0$.

For every $\epsilon\neq\emptyset$ there exists $i\in\{1,\dots,d\}$ with 
$\epsilon_i=1$ and the projection $x\mapsto x_\epsilon$ is invariant 
under $T^*_i$ and thus is $\CW^*$-measurable. 
Therefore we have $\CX^\sharp\subset\CW^*$. We get that 
$\E_{\mu^*}(F\mid\CX^\sharp)=0$. 

 Therefore, by Lemma~\ref{lem:span0} we can restrict to the case that
$$
 F(x)=\prod_{\epsilon\in\{0,1\}^d}
f_\epsilon(x_\epsilon) \text{ where }|f_\epsilon|\leq 1\text{ 
for all }\epsilon\text{ and }\E_\mu(f_{\emptyset}\mid\CZ)=0\ .
$$
We have that $\nnorm{f_\emptyset}=0$ by Lemma~\ref{lem:Z} and 
 by Lemma~\ref{lem:normstar} we have that $\nnorm F^*=0$.\qed

\section{Changing the order of the transformations}
\label{sec:permut}

The next proposition means that we can exchange the order of the limits 
in the formula~\eqref{eq:limit} of Lemma~\ref{lem:limit}.
This result  is parallel to Proposition~3.7 of~\cite{HK1}, but we 
can not simply copy its proof which 
depends of Formula~(9) of~\cite{HK1} which has no analogue in the 
present context. It seems that here we need some 
technology, for example the ``modules'' of~\cite{CL} and/or~\cite{FW}. 
This is the 
only point in this paper where we need  
more elaborate tools.

 \begin{proposition}
\label{prop:permut}
Let $\sigma$ be a permutation of $[d]$, $\sigma_*$ the permutation of  
$\{0,1\}^d$ given by $\bigl(\sigma_*(\epsilon)\bigr)_i=\epsilon_{\sigma(i)}$ for 
every $i$ and  $\Sigma$ 
the associated permutation of $X^*$, given by
$$
 \bigl(\Sigma x\bigr)_\epsilon=x_{\sigma_*(\epsilon)}\text{ for every }
\epsilon\in\{0,1\}^d\ .
$$
Then the box measure associated to the transformations $T_{\sigma(1)},
T_{\sigma(2)},\dots,$ $T_{\sigma(d)}$ is the 
image under $\Sigma$ of the box measure associated to the 
transformations $T_1,T_2,\dots,T_d$.
\end{proposition}

We immediately deduce:
\begin{corollary}
\label{cor:permut}
The seminorm $\nnorm \cdot_{T_1,\dots,T_d}$ remains unchanged if  the 
transformations $T_1,\dots,T_d$ are permuted.  
\end{corollary}
\begin{remark}
A weak form of this Corollary follows easily from Lemma~\ref{lem:limit} 
and Lemma~\ref{lem:uniform} and thus does not depend of the more 
difficult Proposition~\ref{prop:permut}:  The family of functions $f$ 
such that $\nnorm f_{T_1,\dots,T_d}= 0$ does not depend on the order of 
the transformations.
\end{remark}

\subsection{Proof of Proposition~\ref{prop:permut}, first step}
\label{sec:permute}

First we check that it suffices to prove the result for the case of $2$ transformations.
 
Indeed, any permutation of $\{1,\dots,d\}$ can be written as the 
 product of the transposition of two consecutive terms 
 and we can thus assume that $\sigma$ is the 
transposition of  $i$ and $i+1$ for some $i$ with $1\leq i<d$.

Fix $i$ and  let $\tau$ be the box measure associated to $T_1,\dots,T_{i-1}$ (or 
equal to $\mu$ if $i=1$), $S_1=T_i\times\dots\times T_i$ and 
$S_2=T_{i+1}\times\dots\times T_{i+1}$. Applying the result for these 
transformations we get that the box measure associated to 
$T_1,\dots,T_{i-1},T_{i+1},T_i$ is equal to the image of the box 
measure associated to $T_1,\dots,T_{i-1},T_i,T_{i+1}$ under the 
permutation of the last to two digits. We immediately deduce the 
announced result.

Henceforth we assume that $d=2$.
We write $\mu_2$ for the box measure associated to $T_1$ and $T_2$ and 
$\mu_2'$ for the measure associated to $T_2$ and $T_1$ and we want to show 
that 
$\mu'_2$ is the image of $\mu_2$ under the map
$$
 (x_{00},x_{01},x_{10},x_{11})\mapsto 
(x_{00},x_{10},x_{01},x_{11})\colon X^4\to X^4\ .
$$
We recall that 
\begin{gather}
\label{eq:defnu}
 \mu^*=\bigl(\mu\times_{\CI(T_1)}\mu\bigr)
\times_{\CI(T_2\times T_2)}
\bigl(\mu\times_{\CI(T_1)}\mu\bigr)\\
\label{eq:defnup}
 \mu^\circ=\bigl(\mu\times_{\CI(T_2)}\mu\bigr)
\times_{\CI(T_1\times T_1)}
\bigl(\mu\times_{\CI(T_2)}\mu\bigr)
\end{gather}
\subsection{Reduction to the ergodic case}
 We check that we can restrict to the 
case that $(X,\mu,T_1,T_2)$ is ergodic. Indeed, let $\CJ$ be the 
$\sigma$-algebra of sets invariant under $T_1$ and $T_2$ and let
$$
 \mu=\int \mu_\omega\,dP(\omega)
$$
be  the ergodic decomposition of $\mu$ under the action of $T_1$ and 
$T_2$. Since $\CJ\subset\CI(T_1)$ we have that 
$$
 \mu\times_{\CI(T_1)}\mu=\int  \mu_\omega\times_{\CI(T_1)}\mu_\omega
\,dP(\omega)\ .
$$
Since $\CJ\otimes\CJ\subset\CI(T_2\times T_2)$ we have by definition of 
$\mu^*$:
$$
 \mu^*=\int 
\bigl(\mu_\omega\times_{\CI(T_1)}\mu_\omega\bigr)
\times_{\CI(T_2\times T_2)}
\bigl(\mu_\omega\times_{\CI(T_1)}\mu_\omega\bigr)\,dP(\omega)
$$
 and a similar expression holds for $\mu^\circ$. Applying the result to the 
ergodic measures $\mu_\omega$ we deduce the general case.

Henceforth we assume that $(X,\mu,T_1,T_2)$ is ergodic. 
\subsection{Decomposition}
Let $f_{00},f_{10},f_{01}, f_{11}\in L^\infty(\mu)$. We want to show 
that
\begin{multline}
\label{eq:measuresequal}
\int 
f_{00}(x_{00})f_{10}(x_{10})f_{01}(x_{01})f_{11}(x_{11})
\,d\mu^*(x_{00},x_{01},x_{10},x_{11})\\ =
\int f_{00}(x_{00})f_{10}(x_{01})f_{01}(x_{10})f_{11}(x_{11})
\,d\mu^\circ(x_{00},x_{01},x_{10},x_{11})\ .
\end{multline}

Let $\CY$ be the $\sigma$-algebra on $X$ 
corresponding to the maximal isometric extension of 
$(X,\CI(T_1),\mu, T_2)$ in $(X,\mu,T_2)$ and let $\CY'$ be the 
$\sigma$-algebra on $X$ 
corresponding to the maximal isometric extension of 
$(X,\CI(T_2),\mu, T_1)$ in $(X,\mu,T_1)$\footnote{In fact these two 
$\sigma$-algebras are equal but we do not prove this equality here.}.

For every $\epsilon>0$ we can write $f_{00}$ as a sum 
$f_{00}=f+f'+g+h$ of $4$ bounded functions where
$f$ is measurable with respect to $\CY$, $f'$ is measurable with 
respect to $\CY'$, $\E_\mu(g\mid\CY)=\E_\mu(g\mid\CY')=0$ and $\norm 
h_2<\epsilon$. Therefore, we are  reduced to considering  three 
different cases: the case that $f_{00}$ is measurable with respect to 
$\CY$, the completely similar case that $f_{00}$ is measurable with 
respect to $\CY'$, and the case that  $\E_\mu(f_{00}\mid\CY)=
\E_\mu(f_{00}\mid\CY')=0$.
\subsection{The case that $f_{00}$ is measurable with respect to $\CY$}

\begin{lemma}
    \label{lem:uniform2}
Assume that $f_{00}$ is measurable with respect to $\CY$.
 Then 
$$
\sup_{m\in\Z}\Bigl\Vert
\frac 1N\sum_{n=0}^{N-1}
 T_1^n(T_2^mf_{00}\cdot f_{01}) -
\E_\mu\bigl( T_2^m f_{00}\cdot 
f_{01}\mid\CI(T_1)\bigr)\Bigr\Vert_{L^2(\mu)}\to 0
$$
as $N\to+\infty$.
\end{lemma}

\begin{proof}
We use the vocabulary of "modules"  as in~\cite{CL}.
We can restrict to the case that $f_{00}=\phi_i$ where 
$(\phi_1,\dots,\phi_k)$ is a base of a $\bigl(\CI(T_1),T_2\bigr)$-module and 
$1\leq i \leq k$: there exists a $\CI(T_1)$-measurable map $x\mapsto 
U(x)$ with values in the group of unitary $k\times k$ matrices  such 
that 
$$
 T_2\phi_i=\sum_{j=1}^k U_{i,j}\cdot \phi_j \ .
$$
For every $m$,
$$
\E_\mu\bigl( T_2^mf_{00}\cdot f_{01}\mid\CI(T_1)\bigr)
=\sum_{j=1}^kU_{i,j}^{(m)}\cdot 
\E_\mu\bigl(\phi_j f_{01}\mid\CI(T_1)\bigr) 
$$
where $U^{(m)}$ denotes the iterated cocycle:
$$
 U^{(m)}(x)=U(T_2^{m-1}x)\dots U(T_2x)U(x)\ .
$$
For every $n$ 
$$
T_1^n(T_2^m f_{00}\cdot f_{01})= \sum_{j=1}^k U_{i,j}^{(m)}\cdot 
T_1^n(\phi_jf_{01})\ .
$$
 Thus for every $N$ we have
\begin{multline*}
 \Bigl\Vert
\frac 1N\sum_{n=0}^{N-1}  T_1^n(T_2^mf_{00}\cdot f_{01}) -
\E_\mu\bigl( T_2^mf_{00}\cdot f_{01}\mid\CI(T_1)\bigr)
\Bigr\Vert_{L^2(\mu)}\\
\leq \sum_{j=1}^k \Bigl\Vert 
\frac 1N\sum_{n=0}^{N-1}
T_2^n(\phi_jf_{01}) - \E_\mu\bigl(  \phi_j\cdot f_{01}\mid\CI(T_1)\bigr)
\Bigr\Vert_{L^2(\mu)}\ .\qed
\end{multline*}
\renewcommand{\qed}{}
\end{proof}
We now prove formula~\eqref{eq:measuresequal} in the case that 
$f_{00}$ is measurable with respect to $\CY$. 
By Lemma~\ref{lem:limit},
the left hand side is 
equal to
$$
=\lim_{M\to+\infty}\frac 1M\sum_{m=0}^{M-1}
\int\lim_{N\to+\infty}\frac 1N\sum_{n=0}^{N-1}
 T_1^n\bigl(T_2^mf_{00}\cdot f_{01}\bigr)\cdot
 \bigl(T_2^mf_{10}\otimes f_{11}\bigr)\,d\mu\ .
$$
By Lemma~\ref{lem:uniform2}, the limit as $N\to+\infty$ in this 
expression is uniform in $M$, thus the two limits can be permuted 
and the above expression can be rewritten as
$$
 \lim_{N\to+\infty}\frac 1N\sum_{n=0}^{N-1}
\int\lim_{M\to+\infty}\frac 1M\sum_{m=0}^{M-1}
T_2^m\bigl(T_1^nf_{00}\cdot f_{10}\bigr)\cdot
\bigl(T_1^nf_{01}\cdot f_{11}\bigr)\,d\mu
$$
which is equal to the right hand side of~\eqref{eq:measuresequal}.
\qed

\subsection{The case that $\E_\mu(f_{00}\mid\CY)=\E_\mu(f_{00}\mid\CY')=0$}
It is shown in~\cite{CL} that the $T_2\times T_2$ invariant 
$\sigma$-algebra of $(X\times X,\mu\times_{\CI(T_1)}\mu)$ is included 
in $\CY\otimes\CY$. Since $\E_\mu(f_{00}\mid\CY)=0$,  we have
$ \E_{\mu\times_{\CI(T_1)}\mu}\bigl(f_{00}\otimes f_{10}
\mid\CI(T_2\times T_2)\bigr)=0$
and by the definition~\eqref{eq:defnu} of $\mu^*$, 
the left hand side of~\eqref{eq:measuresequal} is equal to zero. By
the same reasoning, the right hand side is also equal to zero.\qed

\end{document}